\documentclass{amsart}
\usepackage{amsmath,amssymb,amscd}\usepackage[all]{xy}

\newtheorem{theorem}{Theorem}
\newcommand{\bt}{\begin{theorem}}
\newcommand{\et}{\end{theorem}}

\newtheorem{lemma}{Lemma}
\newcommand{\bl}{\begin{lemma}}
\newcommand{\el}{\end{lemma}}

\newtheorem{corollary}{Corollary}
\newcommand{\bc}{\begin{corollary}}
\newcommand{\ec}{\end{corollary}}

\newcommand{\beq}{\begin{equation}}
\newcommand{\eeq}{\end{equation}}
\newcommand{\benum}{\begin{enumerate}}
\newcommand{\eenum}{\end{enumerate}}

\newcommand{\C}{\ensuremath{\mathbf C }}

\newcommand{\N}{\ensuremath{ \mathbf N }}

\newcommand{\R}{\ensuremath{\mathbf R}}

\newcommand{\Z}{\ensuremath{\mathbf Z}}

\newcommand{\mbF}{\ensuremath{\mathbf F}}

\newcommand{\mci}{\ensuremath{ \mathcal I}}

\newcommand{\mfI}{\ensuremath{ \mathfrak I}}

\newcommand{\mfp}{\ensuremath{ \mathfrak p}}
\newcommand{\mfP}{\ensuremath{ \mathfrak P}}

\newcommand{\bq}{\begin{eqnarray*}}
\newcommand{\eq}{\end{eqnarray*}}
\newcommand{\be}{\begin{eqnarray}}
\newcommand{\ee}{\end{eqnarray}}
\newcommand{\ba}{\begin{array}}
\newcommand{\ea}{\end{array}}
\newcommand{\bfr}{\begin{flushright}}
\newcommand{\efr}{\end{flushright}}

\newcommand{\bmat}{\left(\begin{matrix}}
\newcommand{\emat}{\end{matrix}\right)}
\newcommand{\bsmallmat}{\left(\begin{smallmatrix}}
\newcommand{\esmallmat}{\end{smallmatrix}\right)}

\DeclareMathOperator{\qqand}{\qquad\text{and}\qquad}

\date{\today}\author{Melvyn B. Nathanson}
\address{Department of Mathematics\\Lehman College (CUNY)\\Bronx, NY 10468}\email{melvyn.nathanson@lehman.cuny.edu}
\title[Krull dimension]{An elementary proof for the Krull dimension of a polynomial ring}


\begin{document}\maketitle

\begin{abstract}. It is proved  that, for every infinite field ${\mbF}$,  the polynomial ring ${\mbF}[t_1,\ldots, t_n]$ has Krull dimension $n$.  The proof uses only ``high school algebra'' and the rudiments of undergraduate  ``abstract algebra.''  \end{abstract}

\section{Prime ideals and  Krull dimension}

In this paper, a \emph{ring} is a commutative ring with a multiplicative identity.  
A \emph{prime ideal} in a ring $R$ is an ideal $\mfP \neq R$ such that, if $a,b \in R$ and $ab \in \mfP$, then $a \in \mfP$ or $b \in \mfP$.
The ring $R$ is an integral domain if and only if $\{ 0\}$ is a prime ideal in $R$.  
In a unique factorization domain, a nonzero principal ideal \mfP\ is prime if and only if 
$\mfP \neq R$ and \mfP\ is generated by an irreducible element.  
For example, in the ring \Z, an ideal \mfP\ is prime if and only if $\mfP = \{ 0\}$ 
or $\mfP = p\Z$ for some prime number $p$.  

An \emph{ideal chain of length $n$}\index{chain!ideal}\index{ideal chain} in the ring $R$ 
is a strictly increasing sequence of $n+1$ ideals of $R$.  
A \emph{prime ideal chain of length $n$}\index{chain!prime ideal}\index{prime ideal chain} 
in $R$ is a strictly increasing sequence of $n+1$ prime ideals of $R$.  
The  \emph{Krull dimension}\index{Krull dimension} of $R$  
is the supremum of the lengths of  prime ideal chains in $R$.   
Eisenbud~\cite[page 215]{eise95} wrote, 
\begin{quotation}
``Arguably the most fundamental notion in geometry and topology is dimension\ldots.  
[Its] \ldots algebraic analogue plays an equally fundamental role in commutative algebra and algebraic geometry.'' 
\end{quotation}
We shall prove that, for every infinite field ${\mbF}$,  
the polynomial ring ${\mbF}[t_1,\ldots, t_n]$ has Krull dimension $n$.  

Krull dimension has many applications.  In algebraic geometry, 
if $S$ is a nonempty set of polynomials in ${\mbF}[t_1,\ldots, t_n]$, 
then the \emph{ variety}\index{variety} (also called the \emph{algebraic set}) 
$V$ determined by $S$ 
is the set of points  in ${\mbF}^n$ 
that are common zeros of the polynomials in $S$: 
\[
V = V(S) =  \left\{ (x_1 , \ldots, x_n ) \in {\mbF}^n: f(x_1,\ldots, x_n) = 0 \text{ for all } f \in S \right\}.
\]
The \emph{vanishing ideal}\index{vanishing ideal}  ${\mfI}(V)$ 
is the set of polynomials that vanish on the  variety $V$:  
\[
{\mfI}(V) =   \left\{ f \in{\mbF}[t_1,\ldots, t_n] : f(x_1,\ldots, x_n) = 0 \text{ for all }  (x_1 , \ldots, x_n ) \in V \right\}.
\]
The quotient ring 
\[
{\mbF}[V] ={\mbF}[t_1,\ldots, t_n]/{\mfI}(V)
\]
is called the \emph{coordinate ring}\index{coordinate ring} of $V$.
One definition of the dimension of the variety $V$ is the Krull dimension of 
its coordinate ring ${\mbF}[V]$.  
For example, if $S = \{0\} \subseteq \mbF[t_1,\ldots, t_n]$ is  the set whose only element 
is the zero polynomial, then $V = V(  \{0\} ) = \mbF^n$, and the vanishing ideal of $V$ 
is $\mci(V) = \mci(\mbF^n) =  \{0\}$.  
We obtain the coordinate ring 
\[
\mbF[V]= {\mbF}[t_1,\ldots, t_n]/{\mfI}(V) \cong {\mbF}[t_1,\ldots, t_n]
\]
and so the variety $\mbF^n$ has dimension $n$.
Nathanson~\cite{nath16} explicitly computes 
the dimensions of some varieties generated by monomials.

The proof that the Krull dimension of the polynomial ring 
${\mbF}[t_1,\ldots, t_n]$ is $n$ uses only high school algebra 
and the rudiments of undergraduate abstract algebra.  
``High school algebra'' means formal manipulations of polynomials. 
Results from high school algebra are collected and proved 
in Sections~\ref{Krull:section:LowerBound} and~\ref{Krull:section:HighSchool}.  
``Rudimentary abstract algebra'' means results whose  proofs use only 
material found in standard undergraduate algebra texts.
Section~\ref{Krull:section:algebra} contains  results from abstract algebra.

For other proofs, see Atiyah and Macdonald~\cite[Chapter 11]{atiy-macd69}, 
Cox, Little, and O'Shea~\cite[Chapter 9]{cox-litt-oshe07}, 
and Kunz~\cite[Chapter 2]{kunz13}.

\section{A lower bound for the Krull dimension}      \label{Krull:section:LowerBound}

A polynomial in one variable $f = \sum_{i =0}^d c_i t^i \in R[t]$ has degree $d$ 
if $c_d \neq 0$.  
The \emph{leading term}\index{leading term} of $f$ is $c_dt^d$, and 
the \emph{leading coefficient}\index{leading coefficient} of $f$ is $c_d$.
The polynomial $f$ is \emph{monic}\index{monic polynomial} 
if its leading coefficient is 1.

Let $\N_0$ denote the set of nonnegative integers.  
For variables $t_1,\ldots, t_n$ and for the $n$-tuple  $I = (i_1, i_2, \ldots, i_n) \in \N_0^n$, 
we define the monomial  
\[
t^I =  t_1^{i_1}t_2^{i_2}\cdots t_n^{i_n}.
\]
The \emph{degree of the monomial}\index{degree!monomial}  $t^I$ is 
$|I| = i_1+i_2+\cdots + i_n$.  
The \emph{degree of the variable}\index{degree!variable} $t_j$ in 
the monomial $t^I$ is $i_j$.

Let ${\mbF}$ be a field.  
The polynomial ring ${\mbF}[t_1,\ldots, t_n]$ is a vector space over ${\mbF}$. 
A basis for this vector space is the set of monomials 
$\left\{ t^I: I \in \N_0^n \right\}$.  
Every nonzero polynomial $f \in \mbF[t_1,\ldots, t_n]$ 
has a unique representation in the form 
\beq         \label{Krull:UniquePolyRep}
f = \sum_{I \in \mci} c_It^I   
\eeq
where $\mci$ is a nonempty finite subset of $\N_0^n$ 
and  $c_I \in {\mbF}\setminus \{ 0\}$ for every $n$-tuple $I \in \mci$.  
The \emph{degree of the polynomial}\index{degree!polynomial}
 $f$ is $\max( |I|: I \in \mci)$, 
and the degree of the variable $t_j$ in $f$ is $\max( i_j: I \in \mci)$.
The degree of the zero polynomial is undefined.  

We write $S \subseteq T$ is $S$ is a subset of $T$, 
and $S \subset T$ if $S$ is a proper subset of $T$.

\bt                          \label{Krull:theorem:polynomial-n}
Let ${\mbF}$ be a field and let $R = {\mbF}[t_1,\ldots, t_n]$ be the polynomial 
ring in $n$ variables with coefficients in \mbF.  Let $\mfP_0 = \{ 0\}$, and,
for $k = 1,2, \ldots, n$, let $\mfP_k$ be the ideal of $R$ 
generated by $\{t_1,t_2,\ldots, t_k\}$. 
Then
\beq                              \label{Krull:polynomial-n}                            
\mfP_0 \subset \mfP_1 \subset \mfP_2 \subset \cdots \subset \mfP_n
\eeq
is a strictly increasing chain of prime ideals in $R$, and so the Krull dimension 
of $R$ is at least $n$.  
\et

\begin{proof}
The ideal $\mfP_0 = \{ 0\}$ is prime because $R$ is an integral domain.  
For $k \in \{1,\ldots, n\}$,  every monomial in every nonzero polynomial in the ideal 
\[
\mfP_k = \left\{ \sum_{j=1}^k t_j f_j: f_j \in R \text{ for } j  \in \{1,\ldots, n\} \right\}  
\]
is divisible by $t_j$ for some $j \in \{1,\ldots, k\}$.
For $k \in \{0,1,\ldots, n-1\}$, we have $t_{k+1} \in \mfP_{k+1}$ but $t_{k+1} \notin \mfP_k$, 
and so~\eqref{Krull:polynomial-n}  is a strictly increasing 
sequence of ideals.  
Moreover, $1 \notin \mfP_n$, and so $\mfP_n$ is a proper ideal in $R$.  
We shall prove that each ideal $\mfP_k$ is prime.  

Let $k \in \{1,\ldots, n\}$.  
Consider the polynomial $f \in R$ of the form~\eqref{Krull:UniquePolyRep}.  
If 
\[
\mci_1 = \{ (i_1,\ldots, i_n)\in \mci:  i_j \geq 1 \text{ for some } j \in \{1,\ldots, k \} \}
\]
and
\[
\mci_2 = \{ (i_1,\ldots, i_n)\in \mci: i_j = 0 \text{ for all } j \in \{1,\ldots, k\}  \}
\]
then $f= f_1 + f_2$, 
where
\[
f_1 =  \sum_{I \in \mci_1}  c_I t^I \in \mfP_k 
\qqand 
f_2 = \sum_{I \in \mci_2}  c_I t^I \in R.
\]
We have $f \in \mfP_k$ if and only if $f_2 = 0$.

For example, let $n = 2$ and  
\beq       \label{Krull:f-example}
f = t_1^3 + 2t_1^2t_2 + 4t_2^3 
\in R = \mbF[t_1,t_2].
\eeq
Choosing $k = 1$, we obtain $\mci_1 = \{ (3,0), (2,1)\}$ 
and $\mci_2 = \{ (0,3) \}$, and so  $f_1 =  t_1^3 + 2t_1^2t_2 \in \mfP_1$ 
and $f_2 =  4 t_2^3  \in R\setminus  \mfP_1$.   
It follows that $f \notin \mfP_1$.  

Let $f \in \mfP_k$, and let $g$ and $h$ be polynomials in $R$ 
such that $f = gh$.  
We write $g = g_1+g_2$ and $h = h_1 + h_2$, where 
$g_1$ and $h_1$ are polynomials in $\mfP_k$, 
and $g_2$ and $h_2$ are polynomials in $R$ that are sums of monomials 
not divisible by $t_j$ for any $j \in \{1,\ldots, k\}$.
Note that  $g_2 h_2 \in \mfP_k$ if and only if $g_2h_2=0$ 
if and only if $g_2 = 0$ or $h_2 = 0$.  
We have  
\[
f = gh = (g_1+g_2) ( h_1 + h_2) = (g_1h_1+g_1h_2+g_2h_1) + g_2h_2.
\]
Because $\mfP_k$ is an ideal, 
\[
g_1h_1+g_1h_2+g_2h_1\in \mfP_k
\]
and so 
\[
g_2 h_2 = f - (g_1h_1+g_1h_2+g_2h_1) \in \mfP_k.  
\]
It follows that either $g_2 = 0$ and $g \in \mfP_k$,  
or $h_2 = 0$ and $h \in \mfP_k$.  
Therefore, $\mfP_k$ is a prime ideal.  
\end{proof}

The main result of this paper (Theorem~\ref{Krull:theorem:KrullDimension}) 
is that the polynomial ring ${\mbF}[t_1,\ldots, t_n]$ has Krull dimension  $n$.  
Although this implies the maximality of the prime ideal chain~\eqref{Krull:polynomial-n},   
there is a nice direct proof of this result.  

\bt                     \label{Krull:theorem:polynomial-max}
The prime ideal chain~\eqref{Krull:polynomial-n}   is maximal.  
\et

\begin{proof}
The first step is to prove that $\mfP_n$ is a maximal ideal.  
The ideal $\mfP_n$ consists of all polynomials with constant term 0.  
Let $g \in R\setminus \mfP_n$, and let \mfI\ be the ideal of $R$ 
generated by the set $\{ t_1,\ldots, t_n, g \}$.
The polynomial $g$ has constant term $c \neq 0$, 
and $g-c \in  \mfP_n \subseteq \mfI$.  
The equation $c^{-1} (g - (g-c) ) = 1$ implies that $1$ is in \mfI, and so $R = \mfI$.  
Thus, the ideal $\mfP_n$ is maximal.  

We shall prove that, for every integer $k \in \{1,2,\ldots, n\}$, 
if $\mfP'$ is a prime ideal of $R$ such that 
\beq               \label{Krull:Krull-q}
\mfP_{k-1} \subset \mfP' \subseteq \mfP_k
\eeq
then  $\mfP' = \mfP_k$.  
This implies that~\eqref{Krull:polynomial-n}  is a maximal prime ideal chain.  

Let $k \geq 1$, and let \mfP'\ be an ideal of $R$ satisfying~\eqref{Krull:Krull-q}.  
It follows that $\mfP' \setminus \mfP_{k-1} \neq \emptyset$.   
Because $\mfP' \subseteq \mfP_k$, 
every polynomial in  $\mfP' \setminus \mfP_{k-1} $ contains at least one monomial 
of the form 
\beq        \label{Krull:KrullDim}
t_k^{i_k} t_{k+1}^{i_{k+1}} \cdots t_n^{i_n}
\eeq
with $i_k \geq 1$.  Let ${\ell}_k$ be the smallest positive integer $i_k$ 
such that a monomial of the 
form~\eqref{Krull:KrullDim} occurs with a nonzero coefficient in some polynomial 
in $\mfP' \setminus \mfP_{k-1} $.  
There exists a  polynomial $f$ 
in $\mfP' \setminus \mfP_{k-1} $ that contains a monomial of the form
~\eqref{Krull:KrullDim} with $i_k = \ell_k$.   
\beq        \label{Krull:KrullDim-2}
t_k^{\ell_k} t_{k+1}^{i_{k+1}} \cdots t_n^{i_n}
\eeq 
We write $f=f_1 + f_2$, where 
\[
f_1 =  \sum_{\substack{I = (i_1,\ldots, i_n) \in \N_0^k \\ i_j \geq 1 
\text{ for some } j \in \{1,\ldots, k-1\} }}  c_I t^I   
\in \mfP_{k-1} 
\]
and
\[
f_2  = \sum_{\substack{I = (i_1,\ldots, i_n) \in \N_0^k \\ i_j = 0 
\text{ for all } j \in \{1,\ldots, k-1\} 
\\ \text{ and } i_k \geq \ell_k}}  c_I t^I  =  t_k^{{\ell}_k} h 
\]
for some nonzero polynomial $h \in R$.  
Because $f_2$ contains the monomial~\eqref{Krull:KrullDim-2} 
in which the variable $t_k$ occurs with degree 
exactly $\ell_k$, the polynomial $h$ must contain a monomial not divisible 
by $t_i$ for all $i \in \{1,\ldots, k-1, k\}$, and so $h \notin \mfP_k$.  
It follows that $ t_k^{{\ell}_k -1} h \notin \mfP'$.

We have $f_2 = f - f_1  \in \mfP'$ because $f \in \mfP'$ and $f_1 \in \mfP_{k-1} \subseteq \mfP'$.  
We factor $f_2$ as follows: 
\[
f_2 =  t_k^{{\ell}_k} h =  t_k \left( t_k^{{\ell}_k -1} h \right).  
\]
Because $ t_k^{{\ell}_k -1} h \notin \mfP'$ and \mfP'\ is a prime ideal,  
it follows that $t_k \in \mfP'$.
Therefore, $\mfP'$ contains $\{t_1,\ldots, t_{k-1}, t_k \}$, and so $\mfP' = \mfP_k$.  
This completes the proof.  
\end{proof}

\section{Results from high school algebra}   \label{Krull:section:HighSchool}

Let $d$ be a nonnegative integer.  
A polynomial $f = \sum_{I \in \mci} c_It^I    \in \mbF[t_1,\ldots, t_n]$  is 
\emph{homogeneous of degree $d$}\index{homogeneous polynomial}\index{polynomial!homogeneous} 
if all of its monomials have degree $d$.  
Let $\lambda \in \mbF$ and $a = (a_1,\ldots, a_n) \in \mbF^n$.
If $t^I = t_1^{i_1}\cdots t_n^{i_n}$ is a monomial of degree $d$, then 
\[
(\lambda a)^I = (\lambda a_1)^{i_1}\cdots (\lambda a_n)^{i_n} 
=  \lambda^{\sum_{j=1}^n i_j}  a_1^{i_1}\cdots a_n^{i_n} 
= \lambda^d a^I.
\]
If $f = \sum_{I \in \mci} c_It^I$ is homogeneous of degree $d$, 
then 
\beq         \label{Krull:homogeneous}
f(\lambda a_1,\ldots, \lambda a_n) =  \sum_{I \in \mci} c_I (\lambda a)^I
=  \lambda^d \sum_{I \in \mci} c_I a^I =  \lambda^df( a_1,\ldots, a_n).
\eeq
For example, the polynomial $f$ defined by~\eqref{Krull:f-example} 
is homogeneous of degree 3, and
\begin{align*}
f(\lambda a_1, \lambda a_2) 
& = (\lambda a_1)^3 + 2(\lambda a_1)^2 (\lambda a_2) + 4(\lambda a_2)^3 \\
& = \lambda^3 \left(  a_1^3 + 2 a_1^2  a_2 + 4 a_2^3 \right) 
= \lambda^3 f( a_1,  a_2).
\end{align*}

Let $ f = \sum_{I \in \mci} c_It^I$ be a nonzero polynomial of degree $d$, 
and let  
\[
\mci_d = \{ I  = (i_1,\ldots, i_n) \in \mci : |I| = \sum_{j=1}^n i_j = d \}.
\]
Then $\mci_d  \neq \emptyset$ and  
\[
f_d =  \sum_{ I \in \mci_d} c_It^I \in {\mbF}[t_1,\ldots, t_n]
\]
is a homogeneous polynomial  of degree $d$.  

The real numbers \R, the complex numbers \C, and the rational functions with 
real coefficients or with complex coefficients are infinite fields.
An example of a finite field is ${\mbF}_2 = \{ 0, 1\}$, 
with addition defined by $0+0=1+1=0$ and $1+0 = 0+1 = 1$, 
and with multiplication defined by $0\cdot 0 = 0\cdot 1 = 1 \cdot 0 = 0$ 
and $1 \cdot 1 = 1$.

In the ring of polynomials $\mbF_2[t_1,t_2]$ with coefficients in $\mbF_2$, 
the homogeneous polynomial $f(t_1,t_2) = t_1^2 + t_1t_2$ has the property that 
\[
f(a_1,1) = a_1^2 + a_1 = 0
\]
for all $a_1 \in \mbF_2$.  The following Lemma shows that this behavior is impossible for 
polynomials with coefficients in an infinite field.

\bl                 \label{Krull:lemma:NoetherNorm1}
Let ${\mbF}$ be an infinite field, and let $f$ be a nonzero polynomial in 
${\mbF}[t_1,\ldots, t_n]$.  
There exist infinitely many points 
$(a_1,\ldots, a_{n-1},a_n) \in \left( \mbF\setminus \{ 0\} \right)^n$ 
such that 
\[
f(a_1,\ldots, a_{n-1},a_n) \neq 0.
\]
If $f$ is homogeneous, then there exist  infinitely many points 
$(a_1,\ldots, a_{n-1}) \in \left( \mbF\setminus \{ 0\} \right)^{n-1}$ such that 
\[
f(a_1,\ldots, a_{n-1},1) \neq 0.
\]
\el

\begin{proof}
By induction on the number $n$ of variables.  
If $n = 1$, then $f(t_1)$ has only finitely many zeros.  Because the field \mbF\ is infinite, 
there exist infinitely many $a_1\in \mbF\setminus \{ 0\}$ with $f(a_1) \neq 0$.  
If $f$ is homogeneous, then $f(t_1) = c_dt_1^d$ for some $d \in \N_0$ 
and $c_d \in  \mbF\setminus \{ 0\}$, 
and $f(1) = c_d \neq 0$.   

Let $n \geq 2$, and assume that the Lemma holds for polynomials in $n-1$ variables.  
A polynomial $f \in {\mbF}[t_1, \ldots, t_{n-1}, t_n]$ 
can also be represented as a polynomial in the variable $t_n$ with coefficients 
in the polynomial ring ${\mbF}[t_1,\ldots, t_{n-1}]$.  
Thus, there exist polynomials $f_0, f_1,\ldots, f_d \in{\mbF}[t_1,\ldots, t_{n-1}]$ such that 
$f_d \neq 0$ and 
\beq        \label{Krull:NoetherNorm}
f(t_1, \ldots, t_{n-1}, t_n) = \sum_{j=0}^d f_j(t_1,\ldots, t_{n-1})  t_n^j.  
\eeq
By the induction hypothesis, there exist infinitely many points 
$(a_1,\ldots, a_{n-1}) \in \left( \mbF\setminus \{ 0\} \right)^{n-1}$ 
such that 
\beq        \label{Krull:NoetherNorm-2}
f_d(a_1,\ldots, a_{n-1}) \neq 0.
\eeq
By~\eqref{Krull:NoetherNorm} and~\eqref{Krull:NoetherNorm-2}, the polynomial
\[
f(a_1, a_2,\ldots, a_{n-1}, t_n) = \sum_{j=0}^d f_j(a_1,\ldots, a_{n-1}) t_n^j \in{\mbF}[t_n]
\]
has degree $d$.    
Because a nonzero polynomial of degree $d$ has at most $d$ roots in ${\mbF}$, 
and because the field ${\mbF}$ is infinite, 
there exist infinitely many elements  $a_n \in{\mbF} \setminus \{0\}$ such that 
\[
f(a_1, a_2,\ldots, a_{n-1}, a_n) \neq 0.
\]

Let $f \in{\mbF}[t_1,\ldots, t_{n-1}, t_n]$ be a homogeneous polynomial of degree $d$, 
and let $(a_1,\ldots, a_{n-1},a_n) \in \left( \mbF\setminus \{ 0\} \right)^n$ 
satisfy $f(a_1,\ldots, a_{n-1},a_n) \neq 0$.  
Applying the homogeneity identity~\eqref{Krull:homogeneous} with $\lambda = a_n^{-1}$, 
we obtain 
\[
f(a_n^{-1} a_1,\ldots, a_n^{-1} a_{n-1}, 1) =  a_n^{-d} f(a_1,\ldots, a_{n-1},a_n) \neq 0
\]
This completes the proof.  
\end{proof}

We shall apply Lemma~\ref{Krull:lemma:NoetherNorm1} to prove 
an important result (Lemma~\ref{Krull:lemma:NoetherNorm2})
that may appear  ``complicated'' and  ``technical,'' but is essential 
in Section~\ref{Krull:section:UpperBound} in the proof of the fundamental 
theorem about Krull dimension.  

Let $t_1,\ldots, t_n, x_1,\ldots, x_n$ be variables, and let $I = (i_1,\ldots, i_n) \in \N_0^n$.  
We consider the polynomials
\[
(t_j+x_jt_n)^{i_j} \in \mbF[x_1,\ldots, x_n,t_1,\ldots, t_n] 
=  \mbF[x_1,\ldots, x_n,t_1,\ldots, t_{n-1}] [ t_n]
\]
for $j = 1, \ldots, n$.  
As a polynomial in $t_n$, the leading term of $(t_j+x_jt_n)^{i_j}$ is 
$x_j^{i_j} t_n^{i_j}$
and the leading term of 
\[
(t_1+x_1t_n)^{i_1} (t_2+x_2t_n)^{i_2}\cdots (t_{n-1}+x_{n-1}t_n)^{i_{n-1}} (x_nt_n)^{i_n}
\]
is 
\beq             \label{Krull:LeadingTerm}
\prod_{j=1}^n (x_jt_n)^{i_j}
=   \left( \prod_{j=1}^nx_j^{i_j} \right) t_n^{\sum_{j=1}^n i_j} 
= x^I  t_n^{ |I| }.
\eeq

Let $g$ be a nonzero polynomial in ${\mbF}[t_1,\ldots, t_n]$.  
If $g$ has degree $d$ in $t_n$, 
then there exist unique polynomials $g_0,\ldots, g_d \in{\mbF}[t_1,\ldots, t_{n-1}]$ 
with $g_d \neq 0$ such that 
\[
g(t_1,\ldots, t_{n-1}, t_n) 
= \sum_{j =0}^d g_j(t_1,\ldots, t_{n-1}) t_n^j.
\] 
The polynomial $g$ is \emph{monic in $t_n$}\index{monic in $t_n$} 
if $g_d =1$, that is, if 
\[
g(t_1,\ldots, t_{n-1}, t_n) 
=  t_n^d + \sum_{j =0}^{d-1} g_j(t_1,\ldots, t_{n-1}) t_n^j.
\]

\bl                 \label{Krull:lemma:NoetherNorm2}
Let ${\mbF}$ be an infinite field, and let $f = \sum_{I \in \mci} c_It^I  \in{\mbF}[t_1,\ldots, t_n]$ 
be a nonzero  polynomial.
There exist $a_1,\ldots, a_{n-1}, \lambda \in{\mbF}$ with $\lambda \neq 0$ 
and a polynomial $g \in{\mbF}[t_1,\ldots, t_n]$ such that $g$ is monic in $t_n$ and 
\[
g(t_1,\ldots, t_{n-1}, t_n)  
= \lambda^{-1} f(t_1+a_1t_n, t_2+a_2t_n, \ldots, t_{n-1}+a_{n-1}t_n, t_n).  
\] 
\el

\begin{proof}
If $f$ has degree $d$, then   
\[
\mci_d = \{ I  = (i_1,\ldots, i_n) \in \mci : |I| = \sum_{j=1}^n i_j = d \} \neq \emptyset 
\]
and  
\[
f_d =  \sum_{ I \in \mci_d} c_It^I \in{\mbF}[t_1,\ldots, t_n]
\]
is a nonzero homogeneous  polynomial of degree $d$.  

We introduce additional variables $x_1,\ldots, x_n$, and consider  
\[
f_d(t_1+  x_1t_n, t_2+x_2t_n, \ldots, t_{n-1}+x_{n-1}t_n, x_nt_n) 
\]
as a polynomial in $t_n$ with coefficients 
in  the ring ${\mbF}[x_1,\ldots, x_n, t_1,\ldots, t_{n-1}]$.
Applying~\eqref{Krull:LeadingTerm}, we have 
\begin{align*}
f_d(t_1+ & x_1t_n, t_2+x_2t_n, \ldots, t_{n-1}+x_{n-1}t_n, x_nt_n) \\
& =  \sum_{I =  (i_1,\ldots, i_n) \in \mci_d } c_I  \prod_{j=1}^{n-1} (t_j+x_jt_n)^{i_j} \  (x_nt_n)^{i_n} \\
& =   \sum_{ I= (i_1,\ldots, i_n) \in \mci_d } c_I \prod_{j=1}^n (x_j t_n)^{i_j}  
+ \text{ lower order terms in $t_n$} \\
& = \left( \sum_{ I \in \mci_d } c_I x^I \right) t_n^d + \text{ lower order terms in $t_n$} \\
& = f_d(x_1,\ldots, x_{n-1}, x_n) t_n^d + \text{ lower order terms in $t_n$}.  
\end{align*}
By Lemma~\ref{Krull:lemma:NoetherNorm1}, there exist $a_1,\ldots, a_{n-1} \in{\mbF}$ 
such that 
\[
\lambda = f_d(a_1,\ldots, a_{n-1},1) \neq 0. 
\]
It follows that 
\begin{align*}
\lambda^{-1}f_d(t_1+ & a_1t_n, t_2+a_2t_n, \ldots, t_{n-1}+a_{n-1}t_n, t_n) \\
& = \lambda^{-1}f_d(a_1,\ldots, a_{n-1}, 1) t_n^d + \text{ lower order terms in $t_n$} \\
& = t_n^d + \text{ lower order terms in $t_n$}
\end{align*}
and so the polynomial 
\begin{align*}
g(t_1, & \ldots, t_{n-1}, t_n)  \\
& = \lambda^{-1}f(t_1+  a_1t_n, t_2+a_2t_n, \ldots, t_{n-1}+a_{n-1}t_n, t_n) \\
& = \lambda^{-1}  f_d(t_1+  a_1t_n, \ldots, t_{n-1}+a_{n-1}t_n, t_n) + \text{ lower order terms in $t_n$} \\
& = \lambda^{-1}  f_d(a_1,\ldots, a_{n-1}, 1) t_n^d + \text{ lower order terms in $t_n$}  \\
& =  t_n^d + \text{ lower order terms in $t_n$}   
\end{align*}
is monic in $t_n$.  
This completes the proof.  
\end{proof}

\section{Results from undergraduate algebra}\label{Krull:section:algebra} 

To obtain an upper bound for the Krull dimension of the polynomial ring, 
we need to study the image of a chain of ideals in a quotient ring.

\bl                \label{Krull:lemma:PQRideals}
If 
\[
\mfI \subset \mfI' \subset \mfI''
\]
is an ideal chain in the ring $R$, then  
\[
\{ \mfI \} \subset \mfI'/\mfI \subset \mfI'' /\mfI 
\]
is an ideal chain in the quotient ring $R/\mfI $. 
If $\mfI'$ is a prime ideal in $R$, then $\mfI'/\mfI$ is a prime ideal in $R/\mfI$.  
\el

\begin{proof}
In the quotient ring $R/\mfI$, the set containing only the coset $0+\mfI = \mfI$ is $\{ \mfI \} $.
The sets $\mfI'/\mfI = \{a  +\mfI: a \in \mfI' \}$ 
and  $\mfI''/\mfI = \{b+\mfI: b \in \mfI''\}$ are ideals in $R/\mfI$.

We have  $\mfI'/\mfI \neq \{ \mfI \} $ because $\mfI' \neq \mfI$.  
Because $\mfI'$ is a proper subset of $\mfI''$, 
there exists  $b \in \mfI''$ with $b \notin \mfI'$, and $b + \mfI \in \mfI''/\mfI$.  
If $b + \mfI \in \mfI'/\mfI$, then there exists $a \in \mfI'$ such that 
$b+ \mfI = a+\mfI$, and so $b-a \in \mfI \subseteq \mfI'$.  
It follows that $b = (b-a)+a \in \mfI'$, which is absurd.  Therefore, $b + \mfI \notin \mfI'/\mfI$
and $\mfI'/\mfI \neq \mfI''/\mfI$.

Let $\mfI'$ be a prime ideal in $R$.  If $a,b \in R$ and 
\[
(a+ \mfI )(b+ \mfI ) \in \mfI'/\mfI
\]
then there exists $c \in \mfI'$ such that $ab+ \mfI = c + \mfI$ and so $ab = c + x$ 
for some $x \in \mfI \subseteq \mfI'$.  Therefore, $ab \in \mfI'$.
Because $\mfI'$ is a prime ideal, we have $a \in \mfI'$ or $b \in \mfI'$, 
and so $a+ \mfI  \in \mfI'/\mfI$ or $b+ \mfI  \in \mfI'/\mfI$.
This completes the proof.  
\end{proof}

\bt   \label{Krull:theorem:strictPQRideals}
If 
\[
\mfP_1 \subset \mfP_2 \subset \cdots \subset \mfP_m 
\]
is a prime ideal chain  in the ring $R$, 
then    
\beq                       \label{Krull:strictPQRideals}
\{ \mfP_1 \} \subset \mfP_2 /\mfP_1 \subset \cdots \subset \mfP_m /\mfP_1 
\eeq
is a prime ideal chain in the quotient ring $R/\mfP_1$.  
\et

\begin{proof}
This follows immediately from Lemma~\ref{Krull:lemma:PQRideals}.
\end{proof}

We also need some results about integral extensions of a ring.    
The ring $S$ is an \emph{extension ring}\index{extension ring} of $R$ 
and the ring $R$ is a \emph{subing}\index{subring} of $S$
if $R \subseteq S$ 
and the multiplicative identity in $R$ is the multiplicative identity in $S$.  
An element $a \in S$ is \emph{integral} over $R$\index{integral element}
if there is a monic polynomial $f \in R[t]$ of degree $d \geq 1$ 
such that $f(a) = 0$.

The ring $S$ is an \emph{integral extension}\index{integral extension ring} of $R$ 
if $S$ is an extension ring of $R$ and every element of $S$ is integral over $R$.  
For example, the ring of Gaussian integers $\Z[i] = \{a+bi:a,b \in \Z\}$
is an integral extension of \Z\ because $a+bi \in \Z[i] $ is a root 
of the monic quadratic polynomial $t^2 -2at + a^2+b^2 \in \Z[t]$.    
The ring $R$ is an integral extension of itself because every element $a \in R$ 
is a root of the monic linear polynomial $t-a \in R[t]$.

\bl         \label{Krull:lemma:extend-1}
Let $S$ be an extension ring of $R$. 
\benum
\item[(i)]
If $\mfP$ is a prime ideal in $S$, then $\mfP \cap R$
is a prime ideal in $R$.  
\item[(ii)]
 Let $S$ be an integral domain that is an integral  extension of the ring $R$.  
If $\mfI$ is a nonzero ideal in $S$, then $\mfI \cap R$ is a nonzero ideal in $R$. 
\eenum 
\el

\begin{proof}
(i)  If $\mfI$ is an ideal in $S$, 
then $\mfI \cap R$ is an ideal in $R$.
Let $\mfP$ be a prime ideal in $S$, and let $\mfp = \mfP \cap R$.  
If $a,b \in R$ and $ab \in \mfp$, then $ab \in \mfP$ and so $a \in \mfP$ 
or $b \in \mfP$.  It follows that $a \in \mfp$ 
or $b \in \mfp$, and so \mfp\ is a prime ideal in $R$.

(ii)  Let $S$ be an integral domain that is an integral  extension of the ring $R$, 
and let $\mfI$ be a nonzero ideal in $S$.  
Let $a \in \mfI$, $a \neq 0$.  
Because $S$ is integral over $R$, there is a monic polynomial $f$ of minimum 
degree $d$ 
\[
f = t^d + c_{d-1}t^{d-1} + \cdots + c_1t + c_0 \in R[t]
\]
such that 
\[
f(a) =  a^d + c_{d-1}a^{d-1} + \cdots + c_1a + c_0 = 0. 
\]
If $c_0 = 0$, then 
\[
\left( a^{d-1} + c_{d-1}a^{d-2} + \cdots + c_1\right) a = 0.
\]
Because $S$ is an integral domain and $a \neq 0$, we obtain
\[
a^{d-1} + c_{d-1}a^{d-2} + \cdots + c_1 = 0
\]
and so $a$ is a root of a monic polynomial of degree $d-1$.  
This contradicts the minimality of $d$.  
Therefore, $c_0 \neq 0$.  
Because $a \in \mfI$ and $\mfI$ is an ideal, we have 
\[
a^d + c_{d-1}a^{d-1} + \cdots + c_1a  
= \left( a^{d-1} + c_{d-1}a^{d-2} + \cdots + c_1\right) a \in \mfI
\]
and so
\[
c_0 = -\left( a^d + c_{d-1}a^{d-1} + \cdots + c_1a  \right) \in \mfI \cap R.
\]
Thus, $\mfI \cap R \neq \{ 0\}$. 
This completes the proof.   
\end{proof}

\bl         \label{Krull:lemma:extend-2}
Let $S$ be a ring, and let $R$ be a subring of $S$ such that $S$ is integral over $R$.
If $\mfP$ and $\mfI$ are  ideals in $S$ such that   
 $\mfP \subset \mfI$ and \mfP\ is prime, 
then $\mfP \cap R \neq \mfI \cap R$.  
\el

\begin{proof}
The quotient ring $S/\mfP$ is an integral domain because the ideal \mfP\ is prime.
The ring $R/\mfP$ is a subring of $S/\mfP$, and $(\mfP \cap R)/\mfP = \{\mfP\}$.  
Thus, to prove that $\mfP \cap R \neq \mfI \cap R$, it suffices 
to prove that $(\mfI \cap R)/\mfP \neq \{\mfP\}$. 

We prove first that  $S/\mfP$ is an integral extension of $R/\mfP$.  
Let $a \in S$, and consider the coset $a + \mfP$.
Because $a$ is integral over $R$, there is a monic polynomial 
\[
f = t^d + \sum_{i=0}^{d-1} c_i t^i \in R[t]
\]
such that 
\[
f(a) = a^d + \sum_{i=0}^{d-1} c_i a^i = 0.
\]
Defining the monic polynomial 
\[
\tilde{f}= t^d + \sum_{i=0}^{d-1} (c_i + \mfP) t^i \in (R/\mfP) [t]
\]
we obtain 
\begin{align*}
\tilde{f}( a + \mfP ) 
& = ( a + \mfP )^d + \sum_{i=0}^{d-1}  (c_i + \mfP)(a+ \mfP)^i  \\
& = \left( a^d + \sum_{i=0}^{d-1} c_i a^i \right) + \mfP \\ 
& = f(a)+\mfP = \mfP
\end{align*}
and so $a + \mfP$ is integral over the ring $R/\mfP$.  
Thus,  $S/\mfP$ is an integral extension of $R/\mfP$.   
 
 The ideal $\mfI/\mfP$ is nonzero in $S/\mfP$ because $\mfP$ is a proper subset of \mfI.
It follows from Lemma~\ref{Krull:lemma:extend-1} that $( \mfI/\mfP)  \cap (R/ \mfP) $
is a nonzero ideal in $R/\mfP$, that is, $( \mfI/\mfP)  \cap (R/ \mfP) \neq \{\mfP\}$. 

We shall prove that $(\mfI \cap R)/\mfP = ( \mfI/\mfP)  \cap (R/ \mfP) $. 
The inclusion  $(\mfI \cap R)/\mfP \subseteq( \mfI/\mfP)  \cap (R/ \mfP) $ 
is immediate.  To prove the opposite inclusion, 
let $r+\mfP \in ( \mfI/\mfP)  \cap (R/ \mfP) $ for some $r \in R$.  
There exists $a \in \mfI$ such that $r+\mfP = a+\mfP$, and so 
 $r -a = b \in \mfP \subset \mfI$.
It follows that $r = a+ b \in \mfI$ and so $r \in \mfI \cap R$.
Therefore, $( \mfI/\mfP)  \cap (R/ \mfP)  \subseteq  (\mfI \cap R)/\mfP$.

We conclude that $( \mfI \cap R) /\mfP = ( \mfI/\mfP)  \cap (R/ \mfP)  \neq \{\mfP\}$,   
and so $\mfP \cap R \neq \mfI \cap R$.
This completes the proof.  
\end{proof}

\bt        \label{Krull:theorem:extend-4}
Let $S$ be an integral domain, and let $R$ be a subring of $S$ 
such that $S$ is integral over $R$.  
If 
\[
\mfP_0 \subset \mfP_1 \subset \cdots \subset \mfP_m
\]
is a prime ideal chain in $S$, then 
\beq          \label{Krull:extend-4}
\mfP_0 \cap R \ \subset \  \mfP_1  \cap R \  \subset \cdots \subset  \ \mfP_m  \cap R 
\eeq
is a prime ideal chain in $R$.  
\et

\begin{proof}
This follows immediately from Lemmas~\ref{Krull:lemma:extend-1} 
and~\ref{Krull:lemma:extend-2}.  
\end{proof}

A \emph{minimal prime ideal}\index{minimal prime ideal}\index{prime ideal!minimal} 
in a  ring $R$ is a nonzero prime ideal \mfP\ such that, if $\mfP' $ is a prime ideal and 
$\{ 0 \} \subseteq \mfP'  \subseteq \mfP$, then $\mfP' = \{ 0 \}$ 
or $\mfP' = \mfP$. 
In a field, the only prime ideal is $\{0\}$ and there is no minimal prime ideal.

\bl                      \label{Krull:lemma:UFD-minimalPrime}
In a unique factorization domain $R$, a prime ideal is minimal if and only if 
it is a principal ideal generated by an irreducible element.  
\el

\begin{proof}
In a unique factorization domain, the principal ideal generated 
by an irreducible element is a nonzero prime ideal.  

Let $\mfP$ be a minimal prime ideal in  $R$.  
Because $\mfP \neq \{0\}$ and $\mfP \neq R$, the ideal $\mfP$ contains 
a nonzero element  that is not a unit.  
This element is a product of irreducible elements.  
Because $\mfP$ is a prime ideal, it contains at least one of these irreducible factors.  
If $a \in \mfP$ and $a$ is irreducible, then $\mfP$ contains the principal ideal  
$\langle a \rangle$, which is a nonzero prime ideal.  
The minimality of $\mfP$ implies that $\mfP = \langle a \rangle$.
Thus, every minimal prime ideal is a principal ideal generated by an irreducible.  

Conversely, let  $a$ be an irreducible element in $R$, and consider the nonzero prime ideal 
$\langle a \rangle$.   If $\mfP$ is a nonzero prime ideal contained in $\langle a \rangle$,   
then $\mfP$ contains a nonzero element that is not a unit.  
This element is a product of irreducibles, and, because 
the ideal $\mfP$ is prime, it must contain at least one of these irreducible elements.  
Let $b$ be an irreducible element in \mfP, 
and let $\langle b \rangle$ be the principal ideal generated by $b$.  
We have
\[
\langle b \rangle \subseteq \mfP \subseteq \langle a \rangle  
\]
and so $a$ divides $b$.  
Because $a$ and $b$ are irreducible elements in $R$, 
it follows that $a$ and $b$ are associates and so 
$\langle a \rangle = \mfP = \langle b \rangle$.  Therefore, 
$ \langle a \rangle$ is a minimal prime ideal.  
This completes the proof.  
\end{proof}

\bl                              \label{Krull:lemma:determinant}
Let $S_0$ be an extension ring of $R$.  
If there is a finite set $\{x_1,\ldots, x_d\} \subseteq S_0$ 
such that every element of $S_0$ is an $R$-linear combination of elements 
of $\{x_1,\ldots, x_d\} $, that is, if
\[
S_0 = \left\{ \sum_{j=1}^d r_j x_j:r_j \in R \text{ for } j = 1,\ldots, d \right\}
\]
then $S_0$  is an integral extension of $R$.  
\el

\begin{proof}
Note that $S_0 \neq \{ 0\}$ because $1 \in S_0$, and so $ \{x_1,\ldots, x_d\} \neq \{0\}$.  

The \emph{Kronecker delta}\index{Kronecker delta} $\delta_{i,j} $ is defined by 
$\delta_{i,j} = 1$ if $i=j$, and $\delta_{i,j} = 0$ if $i \neq j$.  

Let $s \in S_0$.  
Because $S_0$ is a ring and $ \{x_1,\ldots, x_d\} \subseteq S_0$, 
for all $i \in \{1,\ldots, d \}$ we have $sx_i \in S_0$, and so there exist $r_{i,j} \in R$ 
such that 
\[
sx_i =  \sum_{j=1}^d r_{i,j} x_j.
\]
Equivalently, 
\[
\sum_{j=1}^d (\delta_{i,j}s - r_{i,j} )x_j = 0  
\]
and so the homogeneous system of linear equations  
\[
\sum_{j=1}^d (\delta_{i,j}s - r_{i,j} )t_j = 0 \qquad\text{for $i=1,\ldots, d$.}
\]
has the nonzero solution $ \{x_1,\ldots, x_d\}$.  
This implies that the determinant of 
the matrix of coefficients of this system of linear equations is 0.
This $d \times d$ matrix  is 
\[
\bmat 
s - a_{1,1} & -a_{1,2} & \cdots & -a_{1,d}  \\
-a_{2,1} & s - a_{2,2} & \cdots & -a_{2,d}  \\
\vdots & & & \\
-a_{n,1} & -a_{n,2} & \cdots & s - a_{d,d} 
\emat.
\]
and its determinant is a monic polynomial of degree $d$ in $s$ 
with coefficients in $R$.  Therefore, $s$ is integral over $R$. 
This completes the proof.  
\end{proof}

\bl                \label{Krull:lemma:R[a]}
Let $S$ be an integral domain, let $R$ be a subring of $S$, and let $a \in S$.  
Let $R[a]$ be the smallest subring of $S$ that contains $R$ and $a$.
If  $a$ is integral over $R$, then the ring $R[a]$ is an integral extension of $R$.
\el

\begin{proof}
Every element of $R[a]$ is a polynomial in $a$ with coefficients in $R$, 
that is, an $R$-linear combination of elements in the infinite set 
$\{a^i:i=0,1,2,\ldots \}$.  
Because $a$ is integral over $R$, there is a monic polynomial $f \in R[t]$ 
of degree $d$ such that 
$f(a) = 0$.  Rearranging this equation, we obtain 
\[
a^d = \sum_{j=0}^{d-1} c_{d,j} a^j 
\]
with $c_{d,j} \in R$ for $j = 0,1,\ldots, d-1$.  
If $i \geq d$ and 
\[
a^i =  \sum_{j=0}^{d-1} c_{i,j} a^j 
\]
with $c_{i,j} \in R$ for $j = 0,1,\ldots, d-1$, then 
\begin{align*}
a^{i+1} & = a \cdot a^i  = a \sum_{j=0}^{d-1} c_{i,j} a^j 
=  \sum_{j=0}^{d-2} c_{i,j} a^{j+1}+ c_{i,d-1} a^d \\
& =  \sum_{j=1}^{d-1} c_{i,j-1} a^{j}+ c_{i,d-1} \sum_{j=0}^{d-1} c_{i,j} a^j \\
& =   c_{i,d-1} c_{i,0} + \sum_{j=1}^{d-1} ( c_{i,j-1} + c_{i,d-1}  c_{i,j} ) a^j\\
& =   \sum_{j=0}^{d-1} c_{i+1,j} a^j 
\end{align*}
where $c_{i+1,0} =   c_{i,d-1} c_{i,0}$ 
and $c_{i+1,j} =  c_{i,j-1} + c_{i,d-1}  c_{i,j} \in R$ for $j=1,\ldots, d-1$.  
It follows by induction that every nonnegative power of $a$ 
can be written as an $R$-linear combination 
of elements in the finite set $\{1,a,a^2,\ldots, a^{d-1} \}$,  
and so every element in $R[a]$ is also an $R$-linear combination 
of elements in the finite set $\{1,a,a^2,\ldots, a^{d-1} \}$.  
By Lemma~\ref{Krull:lemma:determinant}, the ring $R[a]$ is an integral extension of $R$. 
\end{proof}

\section{An upper bound for the Krull dimension}  \label{Krull:section:UpperBound}

\bt                             \label{Krull:theorem:KrullDimension}
Let ${\mbF}$ be an infinite field.  For every nonnegative integer $n$, 
the Krull dimension of the polynomial ring ${\mbF}[t_1,\ldots, t_n]$  is $n$.  
\et

\begin{proof}
Let $m$ be the Krull dimension of  the ring $R = {\mbF}[t_1,\ldots, t_n]$.  
We have $m \geq n$ by Theorem~\ref{Krull:theorem:polynomial-n}.
We must prove that $m \leq n$.  The proof is by induction on $n$.

If $n = 0$, then $R$ is the field ${\mbF}$, the only prime ideal in a field  is $\{ 0\}$, 
there is no minimal prime ideal, and $R$ has Krull dimension 0.

If $n = 1$, then $R ={\mbF}[t_1]$ is a principal ideal domain, and, consequently, 
a unique factorization domain.   
By Lemma~\ref{Krull:lemma:UFD-minimalPrime}, 
the nonzero prime ideals in $R$ are the principal ideals generated 
by irreducible polynomials,   
and every nonzero prime ideal in $R$  is minimal.  
Therefore, ${\mbF}[t_1]$ has Krull dimension 1.  

Let $n \geq 2$, and let $R =  {\mbF}[t_1,\ldots, t_{n-1}, t_n] = R'[t_n]$, 
where  $R' ={\mbF}[t_1,\ldots, t_{n-1}]$.  
By the induction hypothesis, $R'$ has Krull dimension at most $n-1$.  
Let 
\beq        \label{Krull:MaxChain-1}
\{ 0\} = \mfP_0 \subset \mfP_1 \subset \cdots \subset \mfP_m
\eeq
be a maximal chain of prime ideals in $R$.  
The polynomial ring $R $ is a unique factorization domain, 
and the ideal $\mfP_1$ is a minimal prime ideal in $R$.   
By Lemma~\ref{Krull:lemma:UFD-minimalPrime}, $\mfP_1$ is a principal ideal generated 
by an irreducible polynomial $f = f(t_1,\ldots, t_n)$.  

Here is the critical application of Lemma~\ref{Krull:lemma:NoetherNorm2}:  
There exist $a_1,\ldots, a_{n-1}, \lambda  \in{\mbF}$ with $\lambda \neq 0$ such that 
\[
g = g(t_1,\ldots, t_{n-1}, t_n)  
= \lambda^{-1} f(t_1+a_1t_n, t_2+a_2t_n, \ldots, t_{n-1}+a_{n-1}t_n, t_n) 
\]
is a polynomial that is monic in the variable $t_n$ with coefficients in $R'$.  
We can represent $g$ in the form 
\beq        \label{Krull:tilde-g}
g = \tilde{g}(t_n) = t_n^d + \sum_{i=0}^{d-1} c_i t_n^i \in R'[t_n] 
\eeq
for some positive integer $d$ and polynomials $c_0,c_1,\ldots, c_{d-1} \in R'$.  

The function $\varphi: R \rightarrow R$ defined by $\varphi(a) = a$ for $a \in \mbF$ and 
\[
\varphi(t_j) = \begin{cases}
t_j + a_j t_n & \text{ if $j = 1, \ldots, n-1$} \\
t_n & \text{ if $j=n$} 
\end{cases}
\]
is a ring isomorphism.  Because a ring isomorphism sends prime ideals to prime ideals,  
\[
\{ 0\} = \varphi(\mfP_0) \subset \varphi(\mfP_1) \subset \cdots \subset \varphi( \mfP_m )
\]
is also a maximal chain of prime ideals in $R$.
We have 
\begin{align*}
\varphi(f(t_1,\ldots, t_n) ) & = f(t_1+a_1t_n, t_2 + a_2 t_n, \ldots, t_{n-1} + a_{n-1} t_n, t_n) \\ 
& = \lambda g(t_1,\ldots, t_n)
\end{align*}
and so $\varphi(\mfP_1)$ is the principal ideal generated by $\lambda g(t_1,\ldots, t_n)$.  
Because $\lambda \in \mbF \setminus \{0\}$ is a unit in $R$,  the principal ideal 
$\varphi(\mfP_1)$ is also generated by $g(t_1,\ldots, t_n)$.  
Thus, we can assume that the minimal ideal $\mfP_1$ in the prime ideal 
chain~\eqref{Krull:MaxChain-1}
is a principal ideal generated by a monic polynomial $g \in R'[t_n]$.

The quotient ring $R/\mfP_1$ is an integral domain because $\mfP_1$ is a prime ideal.  
By Theorem~\ref{Krull:theorem:strictPQRideals}, 
\beq        \label{Krull:MaxChain-2}
\{ \mfP_1  \} \subset \mfP_2/\mfP_1  \subset \cdots  \subset \mfP_m/\mfP_1 
\eeq
is a prime ideal chain in  $R/\mfP_1$.
Every coset in the quotient ring $R / \mfP_1$ is of the form 
$f + \mfP_1$, where 
\[
f = \sum_{i=0}^k f_i t_n^i \in R'[t_n]  
\]
and $f_i \in R'$ for $i = 0,1,\ldots, k$.  
It follows that $f_i + \mfP_1 \in R'/\mfP_1$, and so 
\[
f + \mfP_1 = \left( \sum_{i=0}^k f_i t_n^i \right) + \mfP_1 
= \sum_{i=0}^k (f_i + \mfP_1) (t_n + \mfP_1)^i 
\in \left( R'/ \mfP_1 \right) [t_n + \mfP_1].  
\]
This proves that $R / \mfP_1 \subseteq \left( R'/ \mfP_1 \right) [t_n + \mfP_1]$.  Conversely, 
$ \left( R'/ \mfP_1 \right) [t_n + \mfP_1] \subseteq  R / \mfP_1$ and so 
 $\left( R'/ \mfP_1 \right) [t_n + \mfP_1]= R / \mfP_1$.

Thus, we see that the quotient ring $R / \mfP_1$ is also the extension ring of $R' / \mfP_1$ 
that is generated by the coset $t_n + \mfP_1$.  
From~\eqref{Krull:tilde-g}, we have 
\begin{align*}
(t_n + \mfP_1)^d & +  \sum_{i=0}^{d-1} (c_i+ \mfP_1)  (t_n + \mfP_1)^i \\
& = \left(  t_n^d  +  \sum_{i=0}^{d-1} c_i t_n^i \right) + \mfP_1 \\
& = \tilde{g}(t_n) + \mfP_1 = g + \mfP_1= \mfP_1
\end{align*}
and so $t_n + \mfP_1$ is integral over $R' / \mfP_1$.  
By Lemma~\ref{Krull:lemma:R[a]},
$R/\mfP_1 = (R'/\mfP_1) [ t_n + \mfP_1] $  is an integral extension of $R'/\mfP_1$.
By Theorem~\ref{Krull:theorem:extend-4}, 
\beq        \label{Krull:MaxChain-3}
\{ \mfP_1  \}  \cap R'/\mfP_1 \subset \left( \mfP_2/\mfP_1 \right)  \cap R'/\mfP_1 
 \subset \cdots  \subset  \left( \mfP_m/\mfP_1 \right)   \cap R'/\mfP_1 
\eeq
is a  prime ideal chain of length $m-1$ in the ring 
$R'/\mfP_1$.

The degree of $t_n$ in the polynomial $g$ 
is positive, and so the degree of $t_n$ in every nonzero polynomial 
in the principal  ideal $\mfP_1 = \langle g \rangle$ is positive.  
The degree of $t_n$ in every polynomial in $R' = \mbF[t_1,\ldots, t_{n-1}]$ is 0, 
and so
\[
R' \cap \mfP_1 = \{0\}.
\]
This implies that the homomorphism $\psi:R' \rightarrow R/\mfP_1$ 
defined by $\psi(f) =  f+\mfP_1$ is one-to-one, and so 
\[
R' \cong \psi(R') = R'/\mfP_1.
\]
Applying the isomorphism $\psi^{-1}: R'/ \mfP_1 \rightarrow R'$ 
to the maximal prime ideal chain~\eqref{Krull:MaxChain-3} 
in $R'/ \mfP_1$ gives a prime ideal chain of length $m-1$ in $R'$.  
The induction hypothesis implies that $m-1 \leq n-1$. 
This completes the proof.  
\end{proof}


\def\cprime{$'$} \def\cprime{$'$}
\providecommand{\bysame}{\leavevmode\hbox to3em{\hrulefill}\thinspace}
\providecommand{\MR}{\relax\ifhmode\unskip\space\fi MR }
\providecommand{\MRhref}[2]{%
  \href{http://www.ams.org/mathscinet-getitem?mr=#1}{#2}
}
\providecommand{\href}[2]{#2}

\end{document}